\documentclass[11pt]{article}
\usepackage[english]{babel}

\usepackage{color}
\usepackage{array}
\usepackage{enumerate}
\usepackage{enumitem}
\usepackage{refcount}
\usepackage{algorithm}
\usepackage{algorithmic}
\usepackage{boxedminipage}
\usepackage{float}
\usepackage[utf8x]{inputenc}
 
\usepackage{color}
\usepackage{graphicx}
\usepackage{multirow}
\usepackage{tikz}
\usepackage{amssymb}
\usepackage{caption}
\usepackage{enumerate}
\usepackage{enumitem}
\usepackage{amsmath}

\newtheorem{thm}{Theorem}[section]
\newtheorem{prop}[thm]{Proposition}
\newtheorem{lem}[thm]{Lemma}
\newtheorem{cor}[thm]{Corollary}
\newtheorem{defn}[thm]{Definition}

\newtheorem{remark}[thm]{Remark}

\def\proof{{\noindent \bf Proof} }

\newcommand{\bG}{\mathcal G}

\newcommand{\prfend}{\hbox to7pt{\hfil}
\par\vskip-\baselineskip\hbox to\hsize
{\hfil\vbox {\hrule width6pt height6pt}}\vskip\baselineskip}

\DeclareMathOperator{\Spec}{Spec}

\def\a{\bigskip \par \noindent}
\def\b{\par \noindent}
\def\dd{\medskip \par \noindent}
\long\def\eatit#1{}

\def\N{\Bbb N}
\def\C{\Bbb C}
\def\A{\Bbb A}
\def\P{\Bbb P}
\def\J{\Bbb J}
\def\X{\Bbb X}

\def\O{{\cal O}}

\def\m{\mu}

\def\t{\tau}

\font\tengothic=eufm10
\font\sevengothic=eufm7
\newfam\gothicfam
       \textfont\gothicfam=\tengothic
       \scriptfont\gothicfam=\sevengothic
\def\goth#1{{\fam\gothicfam #1}}

\newcommand{\ga}{\goth a}
\newcommand{\gb}{\goth b}
\newcommand{\gm}{\goth m}

\newcommand{\gq}{\goth q}

    \font\tenmsb=msbm10              \font\sevenmsb=msbm7
\newfam\msbfam
       \textfont\msbfam=\tenmsb
       \scriptfont\msbfam=\sevenmsb
\def\Bbb#1{{\fam\msbfam #1}}

\begin{document}
\title{On the Jacobian scheme of a plane curve}
\author{Stefano Canino, Alessandro Gimigliano,
Monica Id\`a\\
}
\maketitle
\abstract{We study the Jacobian scheme of a plane algebraic curve at an ordinary singularity, characterizing it through a geometric property. We compute the Tjurina number for a family of curves at an ordinary singularity showing that it reaches the minimum possible value, using very elementary methods, essentially Gröbner basis. We give an algorithm that gives the analytic type of a double point using the algebraic version of the Mather-Yau Theorem.}
\dd \par \noindent \smallskip {\begin{footnotesize} MSC classification 14B05 14H20, {\it Keywords: Plane curves singularities, Jacobian ideal, Tjurina number.}\end{footnotesize}}
\section{Introduction} 

Let $C: f=0$ be a plane reduced curve passing through a point $P$, let $\goth m_P$ be the maximal ideal of the point and $kP$ the fat point of multiplicity $k$ in $\P^2$, i.e. the 0-dimensional scheme of ideal $\goth m ^k$. Then, to say that $C$ has in $P$  a singular point of multiplicity $m$ means that $C\supset mP$ and $C\not\supseteq (m+1)P$;  this is very rough, but there are others 0-dimensional schemes contained in $C$ which could characterize the singularity more carefully. For example, if $P$ is an $A_n$ singularity, i.e. a double point such that the germ of $C$ at $P$ is analytically equivalent to the germ of the curve $y^2-x^{n+1} =0$ at $O$, then $P$ is a nodal-type singularity (an $A_{n}$ with $n$ odd) if and only if  for any $\ell \geq 1$ there is a curvilinear scheme supported at $P$ of length $\ell$ contained in $C$, while $P$ is a cuspidal singularity $A_{2r}$ if and only if for any $\ell \leq 2r+1$ there is a curvilinear scheme supported at $P$ of length $\ell $ contained in $C$, and there is no curvilinear scheme supported at $P$  of length $> 2r+1$ contained in $C$ (\cite{GI}, Th.2.3). Since the curve $y^2-x^{n+1} =0$ contains the scheme $Y$ of ideal $(y^2, x^{n+1})$, for $n=2r-1$  there is a sort of infinitesimal rectangular pan with an infinite handle contained in the curve, while for $n=2r$ there is no handle, since any curvilinear scheme supported at $P$  of length $ 2r+1$ contained in $C$ is already contained in the rectangular scheme $Y$.

\b So if we want to study a singularity, one of the possible approaches is to answer to the following question: which kind of ``maximal" 0-dimensional scheme supported at $P$ is contained in $C$? But the curve being 1-dimensional, in many cases it will not be possible to bound the length of these schemes, due to the curvilinear schemes contained in $C$. On the other hand, there is a very interesting 0-dimensional subscheme of the curve that gives informations on the singularity, the Jacobian scheme, i.e. the scheme associated to the ideal $(f,f_x,f_y)\C[x,y]_{\goth m_P}$ (see Section 1: note that by some authors, for example in \cite{GLS 1} Definition 2.1 p. 111, this ideal is called the Tjurina ideal, while the Jacobian, or Milnor, ideal for the curve $C:f=0$ at $O$ is defined to be $ (f_x, f_y)\C\{x,y\}$). 
The lengths of the Milnor and of the Jacobian schemes, called respectively the Milnor and  the Tjurina number  and usually denoted by $\m$ and $\t$, have been an intensive object of study in recent years, e.g. see \cite{GLS 2}, \cite{A}, \cite{AABM}, \cite{W}, \cite{HH}, \cite{GH}.
\dd In the case of double points the characterization through the Jacobian scheme is very easy: a double point is of type $A_n$ if if and only if the Jacobian scheme of the curve at $P$ is a curvilinear scheme of length $n$ (see Section 6).  This follows by a theorem of Mather and Yau saying that, if $V: f=0$ and $W: g=0$ are germs of hypersurfaces in $\C^{n+1}$ with isolated singularities at $O$, $(V, O)$ is analytically equivalent to $(W, 0)$ if and only if $ \C\{x_1,\dots,x_{n+1}\}/(f,{\partial f\over \partial x_1},\dots,{\partial f\over \partial x_{n+1}})$ is isomorphic to $ \C\{x_1,\dots,x_{n+1}\}/(g,{\partial g\over \partial x_1},\dots,{\partial g\over \partial x_{n+1}})$ as a $\C$-algebra.
 
 \dd We recall that two germs $(X, p)\subset (\C^{n+1},p)$, $(Y, q) \subset (\C^{n+1},p)$  of hypersurface singularities are said to be topologically, respectively analytically equivalent if there exist two neighborhoods $U$ of $p$ and $V$ of $q$ in $\C^n$, and a map $\phi: U \to V$ mapping $(X, p)$ on $(Y, q)$, with $\phi$ homeomorphism, respectively analytic isomorphism; the corresponding equivalence classes are called topological, respectively analytic types. It is an immediate consequence of the Mather -Yau Theorem that $\t$ is an invariant for analytic types.
\dd The case of double points (Section 6) may lead to think that the length of the Jacobian scheme at a non-ordinary singularity is always bigger than the length at an ordinary singularity (i.e. a singular point of multiplicity $m$ for a curve whose tangent cone at that point consists of $m$ distinct lines), but this is not true. For example, the curves $C: xy(x-y)(x+y)^2+x^6+y^6$ and $D: x^5-y^5=0$ both have multiplicity 5 at $O$, and $O$ is an ordinary multiple point for $D$ while it is non-ordinary for $C$; a computation with CoCoA shows that the length of the Jacobian scheme at $O$ is $15$ for $C$, while for $D$ is clearly 16.

\a In this paper we are interested not only in the length of the Jacobian scheme of a plane curve at a singular point, but expecially in its geometric structure, which is much less studied.
\dd The structure of the paper is as follows: in Section 2 we fix some notations. In Section 3 we show, for lack of references, that if  $\gb:= (g_1,\dots, g_t)\subset  \C[x,y]$ is the ideal of a 0-dimensional scheme with the origin $O$ contained in its support, and $\gm := (x,y)$, then there is a canonical isomorphism of $\C$-algebras
$ \C[x,y]_{\gm}/ \gb  \C[x,y]_{\gm} \cong \C\{x,y\}/ \gb \C\{x,y\}$. This allows us to state the Mather-Yau Theorem \cite{MY} about the classification of isolated hypersurface singularities in a more algebraic setting.

\dd In Section 4 we introduce the notions of $k$-symmetric scheme, i.e. a scheme supported at a point $P$ and intersecting each line  through $P$  with the same length $k$, and of $k$-symmetric local complete intersection ($k$-slci for short) i.e. a local complete intersection of two plane curves having multiplicity $k$ and no tangent in common at $P$. In Theorem \ref{milnor} we show that, if $C$ is a plane curve having an ordinary multiple point at $P$ of multiplicity $m$, the Milnor scheme at $P$ is a  $(m-1)$-slci and the Jacobian scheme is $(m-1)$-symmetric. This gives another way of proving the well known facts that $\mu= (m-1)^2$ and $\tau\leq (m-1)^2$, with the upper bound attained for example by a union of $m$ distinct lines. 
For the equality $\mu= (m-1)^2$, see for example \cite{BK} p.574, where a general formula, which uses the resolution of singularities, is given, or consider that $\m$ is an invariant for topological types, see for example \cite{GLS 1} 3.43.3 p.219, and that any ordinary $m$-multiple point has the same topological type as $x^m+y^m=0$, see for example \cite{GLS 1} 3.30.1 p.202.
\dd For what concerns the lower bound,  in \cite{BGM} the authors compute the minimum of $\t $  in the set of the germs of plane curves topologically equivalent to $x^m+y^n=0$. Keeping in mind that, as just said, any two ordinary $m$-multiple points have the same topological type, and applying \cite{BGM} Proposition p. 550 and Tableau 1 p. 543 for $m=n$, we find that the minimum $\t$ for an ordinary singularity of multiplicity $m$ is $\left\lfloor {3m^2-2m-4\over 4}\right\rfloor $.

\dd The purpose of Section 5 is to use very elementary methods (essentially Gröbner basis) to compute the Tjurina number at $O$ of the plane curves
$$C_{b,c}: x^m+y^m+x^by^c=0, \quad  b+c >m$$  having an ordinary $m$-singularity at $O$, so to prove that (Theorem \ref{poss})
$$ \min_{b,c\, \in \N, b+c>m} \{\tau(C_{b,c})\}= \left\lfloor {3m^2-2m-4\over 4}\right\rfloor.$$
The curves realizing this minimum are, if we ask $b\geq c$, $ x^m+y^m+x^{{m\over 2}+1}y^{m\over 2}=0$ for $m$ even and $ x^m+y^m+x^{m+1\over 2}y^{m+1\over 2}=0$ for $m$ odd. 
\dd In Section 6 we use the Jacobian ideal to give a quick algorithm which determines the type $A_n$ of a double point.
\dd Finally, in the last section we prove a result on the (global) Tjurina number of a curve in $\P^2$; namely, if $C$ is an irreducible curve of degree $d$ and geometric genus $g$, with no infinitely near singular points, then $C$ has only nodes  if and only if $\tau(C)= {d-1 \choose 2}-g$.

\section{Notations and Preliminaries}
In the following the length of a 0-dimensional scheme $Z$ is denoted by  $\ell (Z)$ and
a fat point of multiplicity $m$ in $\P^2$ is denoted by $mP$ (we recall that $\ell(mP)={m+1 \choose 2}$).
Moreover, if $Y$ is a subscheme of $\P^2$, and $P$ is an isolated point of its support,  $Y_P$ denotes the component of $Y$ supported on $P$.
\begin{defn} \rm
\a Let $C: f(x_0,x_1,x_2)=0$ be a reduced curve of degree $d$ in $\P^2$; we denote by $Sing C$ the union of the singular points of $C$. If $P$ is a point of multiplicity $m$ for $C$ , we write $m_P(C)=m$.
\dd In the following we set 
$$\partial _if:= {\partial f \over \partial x_i}, \quad i=0,1,2 $$
and we denote by $C_i$ the degree $d-1$ curves: $$ C_i: \partial _if =0, \quad i=0,1,2 $$
which we call {\em the derivative curves of $C$}.
\a The {\em (projective) Jacobian scheme $\X=\X(C)$ of $C$} is the 0-dimensional subscheme of $\P^2$ defined by the homogeneous, but maybe not saturated, ideal,  called {\em the (projective) Jacobian ideal of $C$}:  
$$\J=\J:=(\partial _0f,\partial _1f,\partial_2f)$$
Hence the support of the Jacobian scheme is $Sing C$, which consists of a finite number of points since the curve is reduced. 
\dd The length $\ell(\X(C))$ of $\X(C)$ is called {\em the global Tjurina number $\tau(C)$ of $C$}.
\end{defn}

\begin{lem}\label{0}  Let $C: f(x_0,x_1,x_2)=0$ be a reduced curve of degree $d$ in $\P^2$ and $P\in C$ with $m_P(C)=m\geq 2$. Then:

\a  a) The curve $C$ contains $\X(C)$.
 \a b)  $m_P(C_i)\geq m-1 $ and for at least one of the $C_i$ the multiplicity at $P$ is exactly $m-1$. 
 \a c) In particular, $\X_P\supseteq (m-1)P$ and $\X_P\not\supseteq mP$.  
\end{lem}
\proof  a) The Euler relation  $d\cdot f=\sum x_i(\partial _if)$ implies that $f\in \J$, that is, $C\supset \X$.
\a b)  $m_P(C)=m \Leftrightarrow {\partial ^{m-1} f \over \partial x_0^jx_1^hx_2^{m-1-j-h}}|_{P}=0$ for $ 0\leq j,h\leq m-1$ and at least one of the derivatives ${\partial ^m f \over \partial x_0^jx_1^hx_2^{m-j-h}}|_{P}$ is $\neq 0$; 
hence for each $i=0,1,2$ and $0\leq j,h\leq m-2$ one has ${\partial ^{m-2} (\partial _i f)\over \partial x_0^jx_1^hx_2^{m-2-j-h}}|_{P}=0$ and for at least one $i$, one $j$, one $h$ we have ${\partial ^{m-1} (\partial _i f)\over \partial x_0^jx_1^hx_2^{m-1-j-h}}|_{P}\neq 0$.
\a c) It is enough to recall that for a curve $D$ we have: $m_P(D)=k \Leftrightarrow D\supseteq kP$ and $D$ does not contain $( k+1)P$. \prfend

\begin{defn} \rm
\a Let $C: f(x,y)=0$ be a reduced curve of degree $d$ in $\A^2$; we use the same notations as in $\P^2$, i.e. $Sing C$ the union of the singular points of $C$, and if $P$ is a point of multiplicity $m$ for $C$, we write $m_P(C)=m$. 
\a In the following we set 
$$f_x:= {\partial f \over \partial x}, \quad f_y:= {\partial f \over \partial y} $$
and we denote by $C_x, C_y$ the degree $d-1$ derivative curves of $C$: $$ C_x: f_x =0, \quad C_y: f_y =0 $$
The {\em (affine) Jacobian scheme $X=X(C)$ of $C$} is the 0-dimensional subscheme of $\A^2$ defined by the ideal, called {\em the (affine) Jacobian ideal of $C$},
$$J=J:=(f,f_x,f_y)$$ 
The curve being reduced, $Sing C$ consists of a finite number of points $P_1,\dots,P_r$, and $X$ is the union of the 0-dimensional schemes $X_{P_1}, \dots, X_{P_r}$. 
\dd {\em The Tjurina number of $C$ at a singular point $P$} is defined to be $\tau(C)_P:=\ell(X_P)$. If no confusion is possible, i.e. when we work locally, looking at the curve $C$ at the point $P$,  we just write $\tau$ instead of $\tau(C)_P$.
\dd It is easy to see that, if $\bar C$ is a curve in $\P^2$, $U_0$ is the affine chart $x_0\neq 0$, and $C=\bar C \cap U_0$, then $\X({\bar C})\cap U_0=X(C)$. Hence, $\X_P=X_P$ for any $P\in Sing C$.
\a The {\em (affine) Milnor scheme $Z=Z(C)$ of $C$} is the subscheme of $\A^2$ defined by the ideal, called {\em the (affine) Milnor ideal of $C$},
$$I=I(C):=(f_x,f_y)$$
Notice that, by \cite{GLS 1} Lemma 2.3 p.113, if $P$ is a singular, necessarily isolated (the curve being reduced), point of $C$, then $P$ is an isolated point of $Z$.
\dd The {\em Milnor number of $C$ at a singular point $P$} is $\m(C)_P:=\ell(Z_P)$. If no confusion is possible, i.e. when we work locally, looking at the curve $C$ at the point $P$,  we just write $\m$ instead of $\m(C)_P$.

\end{defn}

\section{The Mather-Yau Theorem for algebraic curves}
\a In this section we prove, as a consequence of \cite{MY}, that the analytic germs at $O$ of two reduced algebraic plane curves  $C$ and $D$ are biholomorphically equivalent if and only if their algebraic Jacobian schemes $J(C)_O$ and $J(D)_O$ are isomorphic. To do that we need some preliminaries.

\a If $R$ is a ring or a $\C$-algebra, with $\dim R$ we denote the Krull dimension of $R$. If $R$ is a finite $\C$-algebra, i.e. $R$ is finitely generated as $\C$-vector space, the dimension of the $\C$-vector space $R$ is  the length of $R$.
\a We recall that if $Y$ is a closed irreducible subscheme of $\A^n$ or of $\P^n$, one has $\dim Y= \dim \O_{Y,P}$ for any closed point $P\in Y$ (see \cite{Ha} ex.II.3.20 p.94).

\begin{lem}\label{artin}  A noetherian $\C$-algebra $R$ of dimension 0 is a finite $\C$-algebra.
\end{lem}
\proof A ring is an Artin ring if and only if it is noetherian and of dimension 0 (see   \cite{AM} Theorem 8.5 p.90); on the other hand, if $R$ is a finitely generated $k$-algebra ($k$ a field), $R$ is an Artin ring if and only if $R$ is a finite $k$-algebra (see \cite{AM} ex. 8.3 p.92). \prfend

\begin{lem}\label{nota} Let $\gm:=(x,y)$ in $\C[x,y]$, and consider the injective morphisms of $\C$-algebras
$$ \C[x,y] \;\;\buildrel{j}\over \hookrightarrow  \;\;\C[x,y]_{\gm}\;\; \buildrel{\varphi}\over \hookrightarrow  \;\;\C\{x,y\}$$
and set $$ A:=  \C[x,y]_{\gm}, \quad B:= \C\{x,y\}$$
The rings $A$ and $B$ are local rings with maximal ideals respectively $\gm A$ and $\gm B$, and the map $\varphi$ induces an isomorphism on the completions: $\hat \varphi: \hat A \buildrel{\cong}\over \to \hat B$. 
\dd  Moreover, $(A,B)$ is a flat couple, which means that $B$ is a flat $A$-module and for any ideal $\ga$ of $A$ the following holds:
$$\ga B \cap A = \ga$$
\end{lem}
\proof See \cite{GAGA} prop.3 p.9 for the first statement and \cite{GAGA} prop.22 p.36 and prop.28 p.40 for the second one. \prfend

  \begin{thm} Let $g_1,\dots, g_t \in \C[x,y]$, and assume that the subscheme of $\A^2$ associated to the ideal  $\gb:= (g_1,\dots, g_t)$ is 0-dimensional with the origin $O$ contained in its support; let $\gm := (x,y)$.
Then, there is an isomorphism of $\C$-algebras
$$ \C[x,y]_{\gm}/ \gb  \C[x,y]_{\gm} \cong \C\{x,y\}/ \gb \C\{x,y\}$$
 \end{thm}
\proof  
Let $\gb = (g_1,\dots, g_t)= \goth q_1\cap \dots \cap \goth q_n$ be a minimal primary decomposition, with $\goth m_1=\sqrt \goth q_1,  \dots,  \goth m_n= \sqrt \goth q_n$ maximal ideals; since $O$ is in the support, we may assume $\goth m_1=(x,y)$; set $\goth q= \goth q_1$, $\goth m= \goth m_1$. 
\dd Since $\sqrt \goth q=(x,y)$, there exist $n,m > 0$ such that $x^n\in \goth q, y^m\in \goth q$;  hence, every polynomial of degree $\geq n+m-1$ is in $\goth q$, i.e. $\gq \supseteq (x,y)^{n+m-1}$.

\noindent Notice that
$$\gb \C[x,y]_{\gm} = (\goth q\cap \goth q_2\cap  \dots \cap \goth q_n) \C[x,y]_{\gm}= \goth q\C[x,y]_{\gm} \cap \goth q_2\C[x,y]_{\gm} \cap \dots \cap \goth q_n \C[x,y]_{\gm}$$
and $\goth q_j\C[x,y]_{\gm} =\C[x,y]_{\gm}$ for $j=2,\dots,n$, since $\goth q_j\not\subseteq {\gm}$ for $j=2,\dots,n$.
 Hence $\gb \C[x,y]_{\gm} = \gq \C[x,y]_{\gm} $,
so that $$\gb \C\{x,y\} = \gq \C\{x,y\}. \eqno{(\circ)}$$
From now on we use the notations of Lemma \ref{nota}, and we set $\ga:=\gb A$.  
   \b The ideal $\ga B= \gb B$ is the ideal of $B$ generated by $(\varphi \circ j )(g_1),\dots, (\varphi \circ j)(g_t)$, that is, \linebreak  $\ga B = (g_1, \dots, g_t)B= \{ \sigma_1g_1+\dots+\sigma_tg_t, \sigma_1, \dots, \sigma_t \in B  \}$.
 \b Consider the map induced by $\varphi$:
 $$\begin{array}{cccc}
  \Phi: &A/\goth a& \to &B/\goth a B \\
 &f+\goth a& \mapsto &\varphi(f)+ \ga B
 \end{array}$$
 which is well defined, since $f-g \in \ga \Rightarrow \varphi(f)-\varphi(g)= \varphi(f-g) \in \varphi(\ga) \subseteq \ga B$, and is a morphism of $\C$-algebras.
 \b The map $\Phi$ is injective: if  $\varphi(f)\in \ga B$ then, using Lemma \ref{nota}, we have $f\in \varphi ^{-1} (\ga B) = \ga B \cap A=\ga$.
 \a Now we want to prove that $\Phi$ is surjective. 
 \dd Since $A/\ga$ is a noetherian $\C$-algebra of dimension 0, by Lemma \ref{artin}  the $\C$-vector space $A/\ga$ is finite dimensional; let $f_1, \dots, f_s\in A$ such that $ f_1+ \ga, \dots, f_s+ \ga$ is a basis for $A/\ga$. Since $ \C[x,y] \;\;\buildrel{j}\over \hookrightarrow  \;\; A$, the following is true:
 $$ \forall \; q\in  \C[x,y] \;\; \exists\; a_1,\dots,a_s \in \C {\rm \;such \;that 	\;} q-(a_1f_1+ \dots +a_sf_s)\in \ga \eqno{(*)} $$
 If we prove that $ \varphi(f_1)+\ga B, \dots, \varphi(f_s)+\ga B$ is a basis for the $\C$-vector space $B/\ga B$ we are done.  The injectivity of the map $\varphi$ gives the linear independence over $\C$ of $ \varphi(f_1)+\ga B, \dots, \varphi(f_s)+\ga B$.
Now we want to prove that for any  $\sigma \in B$ there are $a_1,\dots, a_s \in \C$ such that
$$ \sigma + \ga B = \sum_{i=1}^s a_i(\varphi(f_i) + \ga B)$$ 
We can write $\sigma = \sum_{0\leq i+j < n+m-1} a_{i,j}x^iy^j  + \sum_{i+j\geq n+m-1} a_{i,j}x^iy^j$; the series $\sum_{i+j\geq n+m-1} a_{i,j}x^iy^j $ being uniformly convergent around $O$ we can rearrange its terms, so we get
$$\sum_{i+j\geq n+m-1} a_{i,j}x^iy^j = x^n( \sum_{i+j\geq n+m-1, i\geq n} a_{i,j}x^{i-n}y^j )+y^m(\sum_{i+j\geq n+m-1, i<n} a_{i,j}x^iy^{j-m}) \;\; \in \gq B$$
By ${(\circ)}$ we have $\ga B =\gb B = \gq B$, so that $$ \sigma + \ga B = \sum_{0\leq i+j < n+m-1} a_{i,j}x^iy^j + \ga B$$ and we conclude by $(*)$ applied to the polynomial $\sum_{0\leq i+j < n+m-1} a_{i,j}x^iy^j$. \prfend

\begin{cor} Let $C: f=0$ be a reduced curve  in $\A^2$ with a singular point in $O$ and Jacobian, respectively Milnor, ideals  $J$, respectively $I$. Let us consider the analytic germ of $C$ at $O$ and let $J^{an}$, respectively $I^{an}$ denote the ideals generated by $f,f_x,f_y$, respectively $f_x,f_y$, in $\C\{x,y\}$. Then there are canonical isomorphisms of $\C$-algebras:

$$ \C[x,y]_{(x,y)}/J\,\C[x,y]_{(x,y)}\cong \C\{x,y\}/J^{an} $$
$$ \C[x,y]_{(x,y)}/I\,\C[x,y]_{(x,y)}\cong \C\{x,y\}/I^{an} $$

\dd In particular, the Tjurina number $\t$, respectively the Milnor number $\m$, of $C$ at $O$ are the dimensions of the analytic algebras $\C\{x,y\} / (f,f_x, f_y)\C\{x,y\}$, respectively $\C\{x,y\} / (f_x, f_y)\C\{x,y\}$.
\end{cor}

 Hence the Theorem in \cite{MY} gives
 
 \begin{thm}\label{M-Y algebrico} Let $C: f=0$, $D:g=0$ be reduced algebraic curves  in $\A^2$ with a singular point at $O$. Then the analytic germs of $C$ and $D$ at $O$ are biholomorphically equivalent if and only if their (algebraic) Jacobian schemes at $O$ are isomorphic as schemes over $\C$.
\end{thm}

\a
\section{Jacobian schemes at ordinary singularities}
 
  In the following, given a polynomial $g$, with $g_k$ we always denote its homogeneous component of degree $k$.
\dd We recall that if  $D$ and $E$ are two plane curves and $P$ is a point  such that $m_D(P)= k, m_E(P)=k$,  $k\geq 1$, then  $(D \cdot E)_P\geq k^2 $, with equality if and only if they do not have common tangents (see for example  \cite{F} property (5) in Section 3.3 or \cite{GLS 1} p.190).
\dd In this section we concentrate our attention on ordinary singularities. Notice that a curve $C$ having a multiple ordinary point at $P$ can have derivative curves with non-ordinary singularities at $P$, as shown in the following remark. 

 \begin{remark} \rm Let $C: f(x,y)=0$ be a reduced curve of degree $d$ in $\A^2$ with $O$ ordinary multiple point, $m_O(C)=m$; it is not restrictive to assume that the line $x=0$ is not a principal tangent at $O$, hence writing $f$ as sum of its homogeneous components we have 
$$f=f_m+ \dots + f_d, \qquad f_m= y^m+\alpha_{m-1}xy^{m-1}+ \dots + \alpha_0x^m$$
and since $f_m$ is the product of $m$ distinct linear factors,  the discriminant $\Delta(g) $ is not zero, where $g(t):= f_m(t,1)$.
Anyhow, the derivative curves $C_x: (f_m)_x+ \dots + (f_d)_x, \quad C_y: (f_m)_y+ \dots + (f_d)_y$ may have a non-ordinary  singularity at $O$, for example
$$f= {1\over 4}y^4-{(2a+b)\over 3}xy^3+ {(2ab+a^2)\over 2}x^2y^2-a^2bx^3y$$ is such that
$$f_y= y^3-(2a+b)xy^2+(2ab+a^2)x^2y-a^2bx^3 = (y-ax)^2(y-bx)$$
hence $C_y$ has a non-ordinary triple point at $O$.
\end{remark}

 \begin{defn}\label{sym} \rm A 0-dimensional scheme $Y$, supported at one point $P\in \P^2$, is said to be  \linebreak {\em $k$-symmetric} when, for every line $r$ passing through $P$, $\ell(Y\cap r)=k$, and
a {\em $k$-symmetric local complete intersection ($k$-slci} for short) if it is a local complete intersection of two curves $D,E$ with no tangent in common at $P$ and such that $m_D(P)= k, m_E(P)=k$, this implying  $\ell (Y) =k^2$.
\end{defn}

\b It is immediate to see that:
\begin{lem}\label{q1bis} A $k$-slci scheme is $k$-symmetric.
\end{lem}

\proof Let $Y$ denote the component of the 0-dimensional scheme of ideal $(\phi, \psi)$ supported on $O$, where 
$$\phi (x,y) = \prod_{j=1}^k  l_j + \phi_{k+1}+ \dots \phi_p, \quad \psi (x,y)= \prod_{j=1}^k  h_j + \psi_{k+1}+ \dots \psi_q .$$
with $l_1=0, \dots, l_k=0, h_1=0,\dots, h_k=0$ lines through $O$, $l_j\neq h_i$ for $i,j=1,\dots,k$. The curves $\phi=0$ and $ \psi=0 $ have no common irreducible component at $O$ (they may obviously have common components away from $O$), since their tangent cones have no common lines. Let $r$ be a line through $O$; choosing coordinates we may assume that  $r: y=0$. We have 
$$\tilde \phi (x):= \phi (x,0) =a_k x^k+a_{k+1}x^{k+1}+\dots + a_sx^s, \tilde \psi (x):= \psi (x,0) = b_k x^k+b_{k+1}x^{k+1}+\dots + b_tx^t$$
where at least one between $a_k$ and $b_k$ is not zero since, by assumption, if $r$ is one of the tangents of the curve $\phi =0$ at $O$, $r$ cannot be in the tangent cone of $\psi =0$. Let's say $a_k\neq 0$; then, there exists a polynomial $f(x)$ with $f(0)\neq 0$ and $\tilde \phi (x)=x^kf(x)$. Moreover, if $\tilde \psi (x) \neq 0$, there exist a polynomial $g(x)$ with $g(0)\neq 0$ and an $n \geq k$ such that $\tilde \psi (x)=x^ng(x)$. Hence
$$ \left (\C[x,y]/(\phi, \psi,\,y)\right)_ {(x,y)} \cong \left (\C[x]/(\tilde \phi, \tilde \psi) \right)_ {(x)} \cong \C[x]/(x^k)$$
\prfend 

 \begin{remark}\label{ksymksqu} \rm $(i)$ A $k$-symmetric scheme needs not to be a $k${\em -slci}, since there are a lot of $k$-symmetric schemes supported on $P$, and the smallest one is $kP$. For example, if $P=O$, all the monomial schemes of ideal $I$ with $$(x^{k}, y^{k}) \subseteq I \subseteq (x,y)^{k}$$ are $k$-symmetric; the above inclusions give  $ {k+1\choose 2} \leq \ell(\C[x,y]/I) \leq k^2$.

  \a $(ii)$ Let us also notice that the scheme $\Spec\left(\C[x,y]/(x^k,y^k)\right)$ is  a $k${\em -slci}, but in general a $k${\em -slci} is not isomorphic to $\Spec\left(\C[x,y]/(x^k,y^k)\right)$. For example, let $C$ and $D$ be the union of the lines  $x=0,y=0,x+y=0, x+2y=0$ and of the lines $x=0,y=0,x+y=0, x+3y=0$ respectively; then, the germs of $C$ and $D$ at $O$ are not analytically equivalent (see \cite{GLS 1} p.157) so that their Jacobian schemes at $O$ are not isomorphic. On the other hand, their Jacobian schemes are $3${\em -slci} by Theorem \ref{milnor} below, hence they cannot be both isomorphic to $Spec\left(\C[x,y]/(x^3,y^3)\right)$.

\a $(iii)$ Let $P$ be a multiple ordinary point of multiplicity $m\leq 3$ for a plane curve $C$; then the Jacobian scheme of $C$ at $P$ is an $(m-1)${\em -slci}. If $m=2$ this follows by Theorem \ref{double}. 
 If $m=3$, by Theorem 2.51 p.152 in \cite{GLS 1},  the germ of $C$ at $P$ is analytically equivalent to the germ of any union of three distinct lines meeting at $O$, for example $D: x^3-y^3=0$; since the Jacobian ideal of $D$ is $(x^2,y^2)$, the conclusion follows by Theorem \ref{M-Y algebrico}.

\a $(iv)$ For any $m$ there exist curves with a multiple ordinary point of multiplicity $m$ at $P$ and such that their Jacobian scheme is a $(m-1)${\em -slci}. In fact, let us consider the curve $C: x^m-y^m=0$, which is the union of the $m$ distinct lines $(x-\eta_jy)=0$ where $\eta_1,\dots,\eta_m$ are the $m^{th} $ roots of unity. Then $C$ has a multiple ordinary point in $O$ of multiplicity $m$, and no other singularities; its Jacobian ideal is
 $$J=(mx^{m-1}, my^{m-1}, x^m-y^m) = (x^{m-1}, y^{m-1})$$
 hence its Jacobian scheme is a $(m-1)${\em -slci}.
 \end{remark}

 \begin{thm}\label{milnor}    Let $P$ be a multiple ordinary point of multiplicity $m$ for a plane curve $C$ and let $Z_P$ be its Milnor scheme at $P$, $X_P$ its Jacobian scheme at $P$. Then:
 \a (i)  the tangent cones of the curves $C_x,C_y$ have no lines in common, hence $Z_P=(C_x \cap C_y)_P$ is a  $(m-1)$-slci;
\dd (ii) $X_P$ is a $(m-1)$-symmetric scheme.
\end{thm}

\proof   $(i)$  It is not restrictive to assume $P=O$. Let $C: f=0$, $f= f_m+\dots +f_{d}$; since $f_m$ is a homogeneous polynomial in $x$ and $y$ of degree $m$, $(f_m)_x$ and $(f_m)_y$ are homogeneous in $x$ and $y$ of degree $m-1$. If $(f_x)_m =0$, we have $f_m=ay^m$ against the assumption that $P$ is an ordinary singularity, and analogously for $f_y$, hence $m_P(C_x)= m_P(C_y)= m-1$.
The polynomials $(f_m)_x$ and $(f_m)_y$ are products of $m-1$ linear factor each, and no factor of $(f_m)_x$ divides $(f_m)_y$  and viceversa. Indeed, assume they do have a linear factor in common: $(f_m)_x=l l_1\dots l_{m-2}, (f_m)_y=l h_1\dots h_{m-2}$; by Euler formula 
$$x(f_m)_x+y(f_m)_y =m\,f_m \; \Rightarrow \; l | f_m  \; \Rightarrow \; f_m=lg  \; \Rightarrow \; (f_m)_x= l_xg + l g_x  \; \Rightarrow \; l | l_xg  \; \Rightarrow \; l |g  \; \Rightarrow \; l^2 | f_m$$
against the assumption that the tangent cone of $C$ at $P$ is the union of $m$ distinct lines.
Hence the curves $C_x$ and $C_y$ have no tangent in common, and we conclude that $Z_P=(C_x \cap C_y)_P$ is a  $(m-1)${\em -slci}.

\a $(ii)$ Let $r:h=0$ be a line passing through $O$;  by $(i)$ and Lemma \ref{q1bis} we have $\ell(Z_P\cap r)=m-1$.  We want to prove that $ \ell(X_P\cap r)= m-1$; we have $(f,f_x,f_y)+(h)=(f,h)+(f_x,f_y,h)$, that is, $X_P\cap r= (C\cap r)\cap (Z_P\cap r)$. Since $P$ has multiplicity $m$ for $C$,  we have $\ell(C\cap r)\geq m$, so that $X_P\cap r$ is obtained intersecting the subscheme of $r$ of length $m-1$ supported on $P$ with a subscheme of $r$ of length $\geq m$ supported on $P$, and the thesis follows.\prfend

\begin{cor}\label{milnor2}   Let $P$ be a multiple ordinary point of multiplicity $m$ for a plane curve $C$. Then $\mu= (m-1)^2$ and $\t \leq (m-1)^2$.
\end{cor} 
\proof  We use the notation of Theorem \ref{milnor}. By Theorem  \ref{milnor} the scheme $Z_P$ is a $(m-1)$-slci, so it has length $\ell(Z_P)= (m-1)^2 $ as observed in Definition \ref{sym}. Moreover, since $Z_P \supseteq X_P$, we have that $\ell(X_P)\leq \ell(Z_P)$.\prfend

\begin{remark}\label{milnor3}  If the plane curve $C: f=0$ is union of $m$ distinct lines through $P$, the polynomial $f$ is homogeneous of degree $m$, so we have $m\,f=xf_x +yf_y$. Hence $(f,f_x,f_y)= (f_x,f_y)$, i.e. the Milnor scheme coincides with the Jacobian scheme and $\t=\m=(m-1)^2$.
\end{remark}

\section{Ordinary singularities with $\tau < \mu$} 
 As recalled in the introduction, the following theorem holds:
\begin{thm}\label{thm}{\rm (\cite{BGM}, \cite{G}, \cite{LP})}
For a plane curve with an ordinary singularity of multiplicity $a$ the Tjurina number takes all the values in the interval $$\left\lfloor {3a^2-2a-4\over 4}\right\rfloor \leq \t  \leq (a-1)^2$$ and the bounds are sharp.
\end{thm}
 The purpose of this section is to show how elementary methods (essentially Gröbner basis) allow us to compute the Tjurina number at $O$ of the plane curves
 $$C_{b,c}: x^a+y^a+x^by^c=0, \quad  b+c >a$$  having an ordinary singularity of multiplicity $a$  at $O=(0,0)$. As a consequence we find that, for any $a\geq 2$,  in this family there are curves which reach the minimum $\t$, i.e.
 $$ \min_{b,c\, \in \N, b+c>m} \{\tau(C_{b,c})\}= \left\lfloor {3a^2-2a-4\over 4}\right\rfloor.$$

\dd Throughout this section $a,b,c$ denote non-negative integers, and we set
$$f=f_{b,c}:=x^a+y^a+x^by^c, \;\quad a\geq 2, \quad  b+c>a, \quad b\geq c;$$ 
clearly, the corresponding results for $b<c$ can be deduced  by symmetry.
\dd  We recall that the ideals $J$ and $\gm$ are defined to be respectively $(f_x,f_y,f)$ and $(x,y)$ in $\C[x,y]$.

\dd For generalities about Gr\"obner basis we refer to  \cite{CLO} Chapter 2, \S 7. With ``grlex order" we mean the graduate lexicographic order. 

\begin{lem}\label{lemmino} We have
	$x^a,y^a\in J$, hence in particular the Jacobian scheme $X= X(C)$ is supported at $O=(0,0)$, and the Jacobian ideal $J$ is primary with radical $\gm$.
\end{lem}
\proof
It is enough to notice that:
$$x^a={1\over (a-b)(a-c)-bc}\left( (a-c)(xf_x-bf)+b(yf_y-cf)\right) $$
$$y^a={1\over (a-b)(a-c)-bc}\left( c(xf_x-bf)+(a-b)(yf_y-cf)\right) $$
\prfend

\begin{prop}\label{taumax} If $b\geq a$ then $J=(x^{a-1},y^{a-1})$ and $\tau=(a-1)^2=\mu$.
\end{prop}
\proof We have 
$$f_x= x^{a-1}(a+bx^{b-a}y^c),\quad f_y= ay^{a-1}+(cx^{b-a+1}y^{c-1})x^{a-1},\quad f=x^{a-1}(x+ x^{b-a+1}y^c)+ y^{a-1}(y)$$
Then, in the local ring $\C[x,y]_\gm$ we have 
 $$ J\C[x,y]_\gm=(x^{a-1}(a+bx^{b-a}y^c), f_y,\, f)\C[x,y]_\gm= (x^{a-1}, f_y,\, f)\C[x,y]_\gm=( x^{a-1}, y^{a-1},\, f)\C[x,y]_\gm= $$ $$=(x^{a-1},y^{a-1})\C[x,y]_\gm$$
 In a minimal primary decomposition of the Jacobian ideal $J$, the primary ideal with radical $\gm$  is the contraction of the ideal $ J\C[x,y]_\gm$ (see \cite{AM}, Proposition 4.8), that is, $J$ being primary, the ideal $(x^{a-1},y^{a-1})$.		
\prfend

\begin{lem}\label{grob}Assume $b<a$ and set
$$f^{(1)}=f_x=bx^{b-1}y^c+ax^{a-1},\quad f^{(2)}=f_y=cx^by^{c-1}+ay^{a-1},\quad f^{(3)}=x^a,\quad f^{(4)}=y^a$$
$$f^{(5)}=x^{a-b-1}y^{a-1},\quad f^{(6)}=x^{a-1}y^{a-c-1},\quad f^{(7)}=x^{a-b}y^{a-1}.$$
Then a reduced Gr\"obner basis of the Jacobian ideal $J$, up to normalization of the leading term and with respect of the grlex order, is given by the following table (the cases marked by "-" cannot occur under our assumptions):
\begin{center}
	{
		\makebox[\textwidth]{{\renewcommand\arraystretch{1.3} 	
		\begin{tabular}{cc|c|c|c|}
			\cline{3-5}
			&   & \textbf{A}                                                            & \textbf{B}                                                            & \textbf{C}                                                        \\ \cline{3-5} 
			&   &  $b={a+1\over 2}$                                                                     & ${a+1\over2}<b<a-1$                                                                     & $b=a-1$                                                                 \\ \hline
			\multicolumn{1}{|c|}{\textbf{1}} & $c<{a-1\over 2}$ & -                                                                     & \begin{tabular}[c]{@{}c@{}}$f^{(1)},f^{(2)},f^{(3)}$\\ $f^{(4)},f^{(7)}$\end{tabular}     & \begin{tabular}[c]{@{}c@{}}$f^{(1)},f^{(2)},f^{(3)}$\\ $f^{(4)},f^{(7)}$\end{tabular} \\ \hline
			\multicolumn{1}{|c|}{\textbf{2}} & $c={a-1\over 2}$ & -                                                                     & \begin{tabular}[c]{@{}c@{}}$f^{(1)},f^{(2)},f^{(3)}$\\ $f^{(4)},f^{(7)}$\end{tabular}     & \begin{tabular}[c]{@{}c@{}}$f^{(1)},f^{(2)},f^{(3)}$\\ $f^{(4)},f^{(7)}$\end{tabular} \\ \hline
			\multicolumn{1}{|c|}{\textbf{3}} & $c={a\over 2}$ & -                                                                     & \begin{tabular}[c]{@{}c@{}}$f^{(1)},f^{(2)},f^{(3)}$\\ $f^{(4)},f^{(5)}$\end{tabular}     & \begin{tabular}[c]{@{}c@{}}$f^{(1)},f^{(3)}$\\ $f^{(5)},f^{(6)}$\end{tabular}     \\ \hline
			\multicolumn{1}{|c|}{\textbf{4}} & $c={a+1\over 2}$ & \begin{tabular}[c]{@{}c@{}}$f^{(1)},f^{(2)},f^{(3)}$\\ $f^{(4)},f^{(5)},f^{(6)}$\end{tabular} & \begin{tabular}[c]{@{}c@{}}$f^{(1)},f^{(2)},f^{(3)}$\\ $f^{(4)},f^{(5)},f^{(6)}$\end{tabular} & \begin{tabular}[c]{@{}c@{}}$f^{(1)},f^{(3)}$\\ $f^{(5)},f^{(6)}$\end{tabular}     \\ \hline
			\multicolumn{1}{|c|}{\textbf{5}} & ${a+1\over 2}<c<a-1$ & - & \begin{tabular}[c]{@{}c@{}}$f^{(1)},f^{(2)},f^{(3)}$\\ $f^{(4)},f^{(5)},f^{(6)}$\end{tabular} & \begin{tabular}[c]{@{}c@{}}$f^{(1)},f^{(3)}$\\ $f^{(5)},f^{(6)}$\end{tabular}     \\ \hline
			\multicolumn{1}{|c|}{\textbf{6}} & $c=a-1$ & - & - & $f^{(5)},f^{(6)}$                                                         \\ \hline
	\end{tabular}}}}
	\nopagebreak\captionof{table}{Gr\"obner basis of $J$.}
	\end{center}
\end{lem} 
\proof Given $g,h\in \C[x,y]$, let $\tilde S(g,h)$ denote the $S$-polynomial of $g$ and $h$ as defined in \cite{CLO}, Chapter 2, \S 6, Definition 4. We set $S(g,h):= \alpha \beta \,\tilde S(g,h)$ where $\alpha$ and $\beta$ are the leading coefficients of $g$ and $h$ (\cite{CLO}, Chapter 2, \S 2, Definition 7). 
\dd
For each case, let us denote by $\bG$ the set of elements appearing in a single cell of Table 1. We denote by $\overline{g}^{\bG}$ the remainder of the division of $g$ by $\bG$ (\cite{CLO}, Chapter 2, \S 3, Theorem 3). It is immediate to see that for monomials $p,q$ one has $\overline{S(p,q)}^\bG=0$. 

\dd We prove that $\bG$ is a reduced Gr\"obner basis using the Buchberger’s Criterion (see \cite{CLO}, Chapter 2, \S 6, Theorem 6), i.e. showing first that $J=(\bG)$, and then that  $\overline{S(f^{(i)},f^{(j)})}^\bG=0\;\forall f^{(i)},f^{(j)}\in\bG$.
We consider the several cases summarized in Table 1, depending on the values of $b$ and $c$.

 \begin{itemize}[leftmargin=*]

\item \textbf{B1) $c<{a-1\over 2}, {a+1\over 2}<b<a-1$}.\\
By Lemma \ref{lemmino}, in order to verify that $J=(\bG)$ is it enough to show that $f\in(\bG)$ and $f^{(7)}\in J$. We have:
$$f={x\over b}f^{(1)}+\left( 1-{a\over b}\right) f^{(3)}+f^{(4)}\Rightarrow f\in(\bG)$$
$$f^{(7)}={1\over a}(x^{a-b}f^{(2)}-cy^{c-1}f^{(3)})\Rightarrow f^{(7)}\in J.$$
Moreover,
$$S(f^{(1)},f^{(2)})=acx^a-aby^a=acf^{(3)}-abf^{(4)},\quad S(f^{(1)},f^{(3)})=ax^{2a-b}=ax^{a-b}f^{(3)}$$
$$S(f^{(1)},f^{(4)})=ax^{a-1}y^{a-c}=y^{a-c}f^{(1)}-bx^{b-1}f^{(4)}$$
$$S(f^{(1)},f^{(7)})=ax^{a-1}y^{a-c-1}={a\over c}x^{a-b-1}y^{a-2c}f^{(2)}-{a^2\over c}x^{a-b-1}y^{a-2c-1}f^{(4)}$$
$$S(f^{(2)},f^{(3)})=ax^{a-b}y^{a-1}=x^{a-b}f^{(2)}-cy^{c-1}f^{(3)},\quad S(f^{(2)},f^{(4)})=ay^{2a-c}=ay^{a-c}f^{(4)}$$
$$S(f^{(2)},f^{(7)})=ay^{2a-c-1}=ay^{a-c-1}f^{(4)},\quad S(f^{(3)},f^{(4)})=S(f^{(3)},f^{(7)})=S(f^{(4)},f^{(7)})=0$$
Thus we have $\overline{S(f^{(i)},f^{(j)})}^\bG=0\;\forall f^{(i)},f^{(j)}\in\bG$, hence $\bG$ is a Groebner basis and it easy to check that it is reduced. 

\item the cases \textbf{C1) $c<{a-1\over 2},b=a-1$}, \textbf{B2) $c={a-1\over 2},{a+1\over 2}<b<a-1$} and \textbf{C2) $c={a-1\over 2},b=a-1$}
are analogous to the case \textbf{B1}.

\item \textbf{B3) $c={a\over 2},{a+1\over 2}<b<a-1$}.\\
 To prove that $J=(\bG)$ we have just to verify that $f^{(5)}\in J$:
$$f^{(5)}=-{c\over a^2}\left( y^{c-1}f^{(1)}-\left( {b\over a}x^{b-1}+{a\over c}x^{a-b-1}\right) f^{(2)}+{bc\over a}x^{2b-a-1}y^{c-1}f^{(3)}\right).$$

As in case \textbf{B1}, one can see that $\overline{S(f^{(1)},f^{(2)})}^\bG=\overline{S(f^{(1)},f^{(3)})}^\bG=\overline{S(f^{(1)},f^{(4)})}^\bG=\overline{S(f^{(2)},f^{(3)})}^\bG=\overline{S(f^{(2)},f^{(4)})}^\bG=0$.  Moreover we have:
$$S(f^{(1)},f^{(5)})=ax^{a-1}y^{a-c-1}=2cx^{a-1}y^{c-1}=2x^{a-b-1}f^{(2)}-2af^{(5)}$$
$$S(f^{(2)},f^{(5)})=ay^{2a-c-1}=ay^{a-c-1}f^{(4)},\quad S(f^{(3)},f^{(4)})=S(f^{(3)},f^{(5)})=S(f^{(4)},f^{(5)})=0$$
so that $\overline{S(f^{(i)},f^{(j)})}^\bG=0,\;\forall f^{(i)},f^{(j)}\in\bG$. Hence $\bG$ is a Gr\"obner basis and again it is easy to check that it is reduced. 

\item \textbf{C3) $c={a\over 2},b=a-1$}.\\
The proof that $(f^{(1)},f^{(2)},f^{(3)},f^{(4)},f^{(5)})$ is a Gr\"obner basis  is analogous to the previous one. However this is not a reduced Gr\"obner basis. Indeed, we have $f^{(5)}=x^{a-b-1}y^{a-1}=y^{a-1}$ so that we can remove $f^{(4)}$ and replace $f^{(2)}$ by $f^{(2)}-af^{(5)}=x^by^{c-1}=x^{a-1}y^{a-c-1}=f^{(6)}$. Hence the reduced Gr\"obner basis is $\bG=(f^{(1)},f^{(3)},f^{(5)},f^{(6)})$.
\item \textbf{B5) ${a+1\over 2}<c<a-1,{a+1\over 2}<b<a-1$}.\\
In order to prove that $J=(\bG)$ it is enough to verify that $f^{(5)},f^{(6)}\in J$. In fact, we have:
$$f^{(5)}={1\over a^2}\left( -cy^{c-1}f^{(1)}+ax^{a-b-1}f^{(2)}+bcx^{b-1}y^{2c-a-1}f^{(4)}\right) $$
$$f^{(6)}={1\over a^2}\left( ay^{a-c-1}f^{(1)}-bx^{b-1}f^{(2)}+bcx^{2b-a-1}y^{c-1}f^{(3)}\right).$$
As in case \textbf{B1}, one can see that $\overline{S(f^{(1)},f^{(2)})}^\bG=\overline{S(f^{(1)},f^{(3)})}^\bG=\overline{S(f^{(1)},f^{(4)})}^\bG=\overline{S(f^{(2)},f^{(3)})}^\bG=\overline{S(f^{(2)},f^{(4)})}^\bG=0$.  Moreover we have:
$$S(f^{(1)},f^{(5)})=ax^{a-1}y^{a-c-1}=af^{(6)},\quad S(f^{(1)},f^{(6)})=ax^{2a-b-1}=ax^{a-b-1}f^{(3)}$$
$$S(f^{(2)},f^{(5)})=ay^{2a-c-1}=ay^{a-c-1}f^{(4)},\quad  S(f^{(2)},f^{(6)})=ax^{a-b-1}y^{a-1}=af^{(5)}$$
$$S(f^{(3)},f^{(4)})=S(f^{(3)},f^{(5)})=S(f^{(3)},f^{(6)})=S(f^{(4)},f^{(5)})=S(f^{(4)},f^{(6)})=S(f^{(5)},f^{(6)})=0.$$
Hence $\bG$ is a Gr\"obner basis and, again, it is easy to check that it is reduced. 

\item the cases \textbf{A4) $c={a+1\over 2},b={a+1\over 2}$} and \textbf{B4) $c={a+1\over 2},{a+1\over 2}<b<a-1$}
are analogous to the case \textbf{B5}.

\item \textbf{C5) ${a+1\over 2}<c<a-1, b=a-1$}.\\
The proof that $(f^{(1)},f^{(2)},f^{(3)},f^{(4)},f^{(5)},f^{(6)})$ is a Gr\"obner basis is analogous to the case \textbf{B5}, but this is not a reduced one. To show that a reduced Gr\"obner basis is $\bG=(f^{(1)},f^{(3)},f^{(5)},f^{(6)})$ one can proceed as in case \textbf{C3}.
\item \textbf{C4) $c={a+1\over 2},b=a-1$}.\\
It is analogous to the case \textbf{C5}.
\item \textbf{C6) $c=a-1,b=a-1$}.\\
This case follow easily by case \textbf{C5} noting that $f^{(5)}=y^{a-1},f^{(6)}=x^{a-1}$.\prfend
\end{itemize}
\pagebreak
\begin{cor}\label{leading}
Assume $b<a$; then a system of generators for the leading terms ideal $(LT(J))$ of $J$ is given by the following table:
\begin{center}
	{
		\makebox[\textwidth]{{\renewcommand\arraystretch{1.5} 	
				\begin{tabular}{cc|c|c|c|}
					\cline{3-5}
					&   & \textbf{A}                                                            & \textbf{B}                                                            & \textbf{C}                                                        \\ \cline{3-5} 
					&   & $b={a+1\over 2}$                                                                     & ${a+1\over2}<b<a-1$                                                                     & $b=a-1$                                                                 \\ \hline
					\multicolumn{1}{|c|}{\textbf{1}} & $c<{a-1\over 2}$ & -                                                                     & \begin{tabular}[c]{@{}c@{}}$x^{b-1}y^c,x^by^{c-1}$\\ $x^a,y^a,x^{a-b}y^{a-1}$\end{tabular}     & \begin{tabular}[c]{@{}c@{}}$x^{b-1}y^c,x^by^{c-1}$\\ $x^a,y^a,x^{a-b}y^{a-1}$\end{tabular} \\ \hline
					\multicolumn{1}{|c|}{\textbf{2}} & $c={a-1\over 2}$ & -                                                                     & \begin{tabular}[c]{@{}c@{}}$x^{b-1}y^c,x^by^{c-1}$\\ $x^a,y^a,x^{a-b}y^{a-1}$\end{tabular}     & \begin{tabular}[c]{@{}c@{}}$x^{b-1}y^c,x^by^{c-1}$\\ $x^a,y^a,x^{a-b}y^{a-1}$\end{tabular} \\ \hline
					\multicolumn{1}{|c|}{\textbf{3}} & $c={a\over 2}$    & -                                                                     & \begin{tabular}[c]{@{}c@{}}$x^{b-1}y^c,x^by^{c-1}$\\ $x^a,y^a,x^{a-b-1}y^{a-1}$\end{tabular}     & \begin{tabular}[c]{@{}c@{}}$x^{b-1}y^c,x^{a-b-1}y^{a-1}$\\ $x^a,x^{a-1}y^{a-c-1}$\end{tabular}     \\ \hline
					\multicolumn{1}{|c|}{\textbf{4}} & $c={a+1\over 2}$ & \begin{tabular}[c]{@{}c@{}}$x^{b-1}y^c,x^{b}y^{c-1},x^a,y^a$\\ $x^{a-b-1}y^{a-1},x^{a-1}y^{a-c-1}$\end{tabular} & \begin{tabular}[c]{@{}c@{}}$x^{b-1}y^c,x^{b}y^{c-1},x^a,y^a$\\ $x^{a-b-1}y^{a-1},x^{a-1}y^{a-c-1}$\end{tabular} & \begin{tabular}[c]{@{}c@{}}$x^{b-1}y^c,x^{a-b-1}y^{a-1}$\\ $x^a,x^{a-1}y^{a-c-1}$\end{tabular}     \\ \hline
					\multicolumn{1}{|c|}{\textbf{5}} & ${a+1\over 2}<c<a-1$ & - & \begin{tabular}[c]{@{}c@{}}$x^{b-1}y^c,x^{b}y^{c-1},x^a,y^a$\\ $x^{a-b-1}y^{a-1},x^{a-1}y^{a-c-1}$\end{tabular} & \begin{tabular}[c]{@{}c@{}}$x^{b-1}y^c,x^{a-b-1}y^{a-1}$\\ $x^a,x^{a-1}y^{a-c-1}$\end{tabular}     \\ \hline
					\multicolumn{1}{|c|}{\textbf{6}} & $c=a-1$ & - & - & \begin{tabular}[c]{@{}c@{}}$x^{a-b-1}y^{a-1}$\\ $x^{a-1}y^{a-c-1}$\end{tabular}                                                         \\ \hline
	\end{tabular}}}}
	\nopagebreak\captionof{table}{Generators for  $(LT(J))$.}
	\end{center}
\end{cor}
\proof It follows by the definition of Gr\"obner basis and by Lemma \ref{grob}. \prfend

\begin{prop}\label{lunghezze}
Assume $b<a$ and set
$$\ell_1(a,b,c)=b(a-1)+c(a+1)-bc-a+1,\quad \ell_2(a,b,c)=b(a-1)+c(a+1)-bc-a$$
$$\ell_3(a,b,c)=(a-1)^2,\quad \ell_4(a,b,c)=b(a-1)+c(a-1)-bc$$
Let $Y$ be the scheme associated to $(LT(J))$; then $\ell(Y)$ is given by the following table:
\begin{center}
	{
		\makebox[\textwidth]{{\renewcommand\arraystretch{1.5} 	
				\begin{tabular}{cc|c|c|c|}
					\cline{3-5}
					&   & \textbf{A}                                                            & \textbf{B}                                                            & \textbf{C}                                                        \\ \cline{3-5} 
					&   & $b={a+1\over 2}$                                                                     & ${a+1\over2}<b<a-1$                                                                     & $b=a-1$                                                                 \\ \hline
					\multicolumn{1}{|c|}{\textbf{1}} & $c<{a-1\over 2}$  & -                                                                     & $\ell_1(a,b,c)$     & $\ell_1(a,b,c)$ \\ \hline
					\multicolumn{1}{|c|}{\textbf{2}} & $c={a-1\over 2}$ & -                                                                     & $\ell_1(a,b,c)$    & $\ell_1(a,b,c)$ \\ \hline
					\multicolumn{1}{|c|}{\textbf{3}} & $c={a\over 2}$ & -                                                                     & $\ell_2(a,b,c)$     & $\ell_3(a,b,c)$     \\ \hline
					\multicolumn{1}{|c|}{\textbf{4}} & $c={a+1\over 2}$ & $\ell_4(a,b,c)$ & 
				$\ell_4(a,b,c)$ & $\ell_3(a,b,c)$     \\ \hline
					\multicolumn{1}{|c|}{\textbf{5}} & ${a+1\over 2}<c<a-1$ & - & $\ell_4(a,b,c)$ & $\ell_3(a,b,c)$     \\ \hline
					\multicolumn{1}{|c|}{\textbf{6}} & $c=a-1$ & - & - & $\ell_3(a,b,c)$                                                         \\ \hline
	\end{tabular}}}}
\nopagebreak\captionof{table}{Length of $Y$.}
	
\end{center}
\end{prop}

\proof The length of a scheme given by a monomial ideal can be easily computed if we know a system of generators for the ideal, for example using its graphic representation (see \cite{CLO}, Chapter 9, \S 2, Example 1). In our case, generators of $(LT(J))$ are given by Corollary \ref{leading}, and it is enough to prove cases \textbf{B1},\textbf{B3},\textbf{B5},\textbf{C5} and \textbf{C6} because the other ones are analogous. The graphic representation of $Y$ in case \textbf{B1}, \textbf{B3}, \textbf{B5} and \textbf{C5} is given in the figure below, starting on the left upper corner and going clockwise:
 
\begin{center}
\scalebox{.6}{\includegraphics{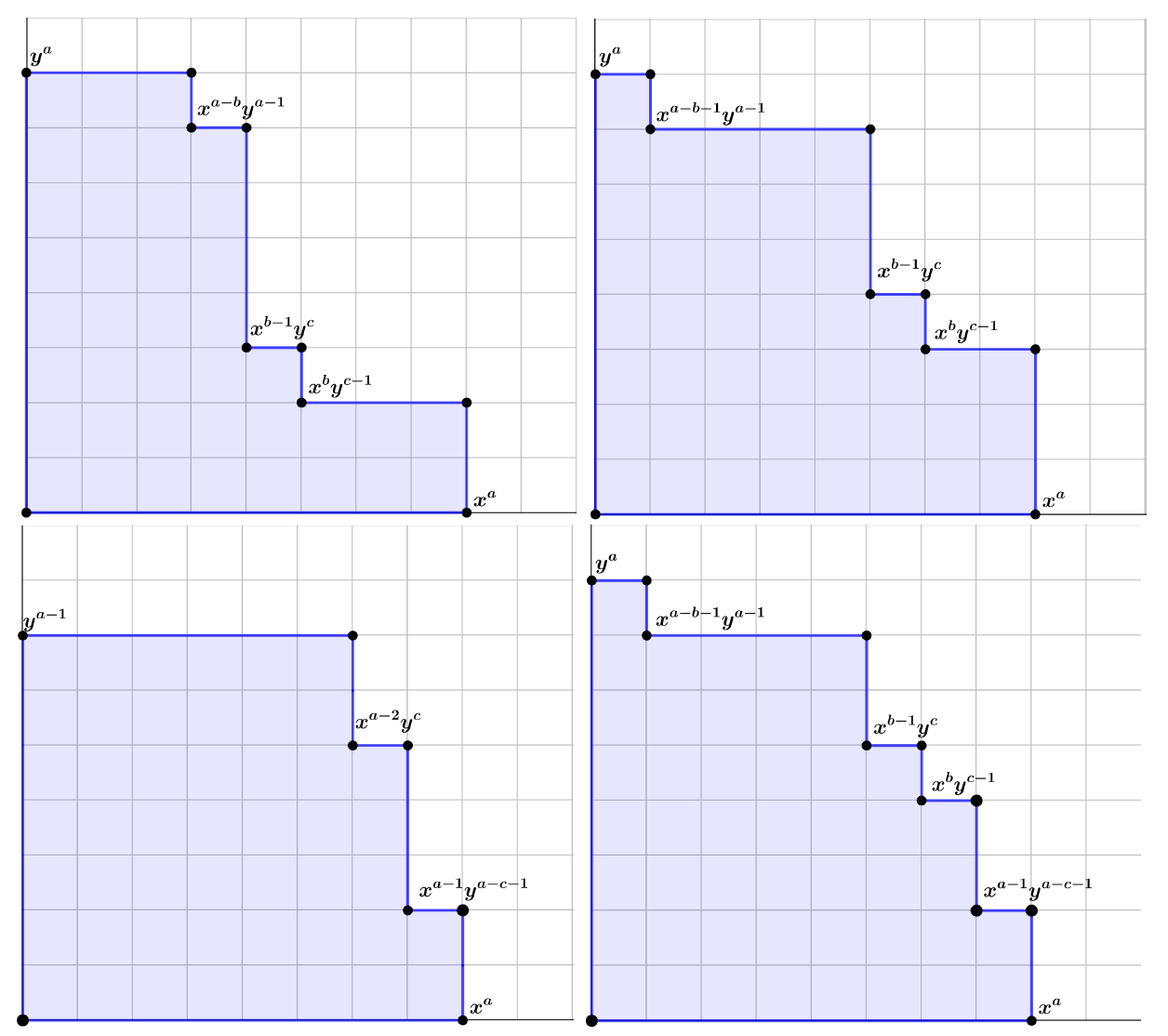}}
\nopagebreak\captionof{figure}{The blue areas give $\ell(Y)$.}
\end{center}

\a Case \textbf{B1 ($c<{a-1\over 2}, {a+1\over 2}<b<a-1$}): since $(LT(J))=(x^{b-1}y^c,x^by^{c-1},x^a,y^a,x^{a-b}y^{a-1})$, 
$$\ell(Y)=(a-b)a+(2b-a-1)(a-1)+c+(a-b)(c-1)=b(a-1)+c(a+1)-bc-a+1.$$
Case \textbf{B3 ($c={a\over 2},{a+1\over 2}<b<a-1$}): since $(LT(J))=(x^{b-1}y^c,x^by^{c-1},x^a,y^a,x^{a-b-1}y^{a-1})$,
$$\ell(Y)=(a-b-1)a+(2b-a)(a-1)+c+(a-b)(c-1)=b(a-1)+c(a+1)-bc-a.$$
Case \textbf{B5 (${a+1\over 2}<c<a-1,{a+1\over 2}<b<a-1$}): since $(LT(J))=(x^{b-1}y^c,x^by^{c-1},x^a,y^a,x^{a-b-1}y^{a-1},x^{a-1}y^{a-c-1})$, 
$$\ell(Y)=(a-b-1)a+(2b-a)(a-1)+c+(a-1-b)(c-1)+a-c-1=b(a-1)+c(a-1)-bc.$$
Case \textbf{C5 (${a+1\over 2}<c<a-1, b=a-1$}): since $(LT(J))=(x^{a-2}y^c,x^a,y^{a-1},x^{a-1}y^{a-c-1})$, 
$$\ell(Y)=(a-2)(a-1)+c+a-c-1=(a-1)^2.$$
Case \textbf{C6 ($c=a-1,b=a-1$}): since $(LT(J))=(x^{a-1},y^{a-1})$, $\ell(Y)= (a-1)^2.$
 \prfend

\begin{thm}\label{poss}
Let $a\in\N$ with $a\geq 2$ and $b,c \in \N$ with $b+c >a$. Then the curve
 $$C_{b,c}:\quad  x^a+y^a+x^by^c=0$$ has an ordinary singularity of multiplicity $a$  at $O$, and, under the assumption $ b\geq c$, the possible values for its Tjurina number $\tau$ at $O$  are: 
$$
\tau= \left\{ \begin{array}{l}
(a-1)^2 \quad  if\; b\geq a \; or \; b=a-1,\; c\geq {a\over 2}\\
\\
 b(a-1)+c(a+1)-bc-a+1\quad  if \; {a+1\over 2} < b \leq a-1,\; c\leq {a-1\over 2}\\
 \\
 b(a-1)+c(a+1)-bc-a \quad  if \; {a+1\over 2} < b < a-1,\; c= {a\over 2}\\
 \\
 b(a-1)+c(a-1)-bc \quad  if \; {a+1\over 2} \leq b < a-1, \; {a+1\over 2} \leq c < a-1\\
\end{array}\right.
$$
\a The minimum value reached by the Tjurina number $\tau$ of $C_{b,c}$ at $O$ is 
 $$\tau_a:= \min_{b,c\, \in \N, b+c>a} \{\tau(C_{b,c})_O\}= \left\lfloor {3a^2-2a-4\over 4}\right\rfloor$$
and the curves realizing this minimum are, if we ask $b\geq c$, $ x^m+y^m+x^{{m\over 2}+1}y^{m\over 2}=0$ for $m$ even and $ x^m+y^m+x^{m+1\over 2}y^{m+1\over 2}=0$ for $m$ odd.
\end{thm}
\a
\proof If $Y$  denotes the scheme associated to the leading terms ideal $(LT(J))$, then the affine Hilbert function of $\C[x,y]/J$ is equal to the affine Hilbert function of $\C[x,y]/(LT(J))$ (see \cite{CLO}, Chapter 9, \S 3, Proposition 4), so that, $Y$ being supported uniquely at $O$, $\tau= \ell(Y)$. Hence the first statement follows from Lemma \ref{taumax} and Proposition \ref{lunghezze}.
\dd For the second statement, we prove, under the assumption $\; b\geq c$, that $\tau_a={3a^2-2a-4\over 4}$ if $a$ is even and $\tau_a= {3a^2-2a-5\over 4}$ if $a$ is odd, and once this is done we can drop the assumption $\; b\geq c$ just exchanging the variables $x$ and \nolinebreak $y$.
\dd If $P$ be a multiple ordinary point of multiplicity $a\leq 3$ for any plane curve $D$, we already noticed in Remark \ref{ksymksqu}  $(iii)$ that the Jacobian scheme of $D$ at $P$ is an $(a-1)${\em -slci}, hence $\tau=1$ if $a=2$, $\tau=4$ if $a=3$, and in fact $\tau_2=1$ and $\tau_3=4$. Hence in the following we may assume $a\geq 4$.

 \dd We use the notations of Proposition \ref{lunghezze}; recall that $\ell_1=b(a-1)+c(a+1)-bc-a+1,\;\ell_2=b(a-1)+c(a+1)-bc-a,\;  \ell_4=b(a-1)+c(a-1)-bc$. 
 \dd Assume $a$ even. Since $(a-1)^2$ is the greatest possible value for an ordinary singularity of multiplicity $a$, it is enough to minimize $\ell_1,\ell_2$ and $\ell_4$ in their definition domains. We find the minimum of $\ell_1,\ell_2$ and $\ell_4$ and finally the minimum between these three minima.
	Let us start with $\ell_1$. Its domain is defined by the following inequalities: 
	$$b+c\geq a+1,\quad  {a+1\over 2}<b\leq a-1,\quad  c\leq {a-1\over 2} \quad ({\rm hence} \;\;   b\geq c)$$ 
but, since $a$ is even and $a,b,c\in\N$, we can refine these inequalieties to the following ones:
	$$b+c\geq a+1,\quad  {a\over 2}+1\leq b\leq a-1,\quad 2\leq c\leq {a\over 2}-1.$$
	The plane region corresponding to these inequalities is the triangle $T$ as shown by Fig. 2:
	\begin{center}
		\includegraphics[scale=0.7]{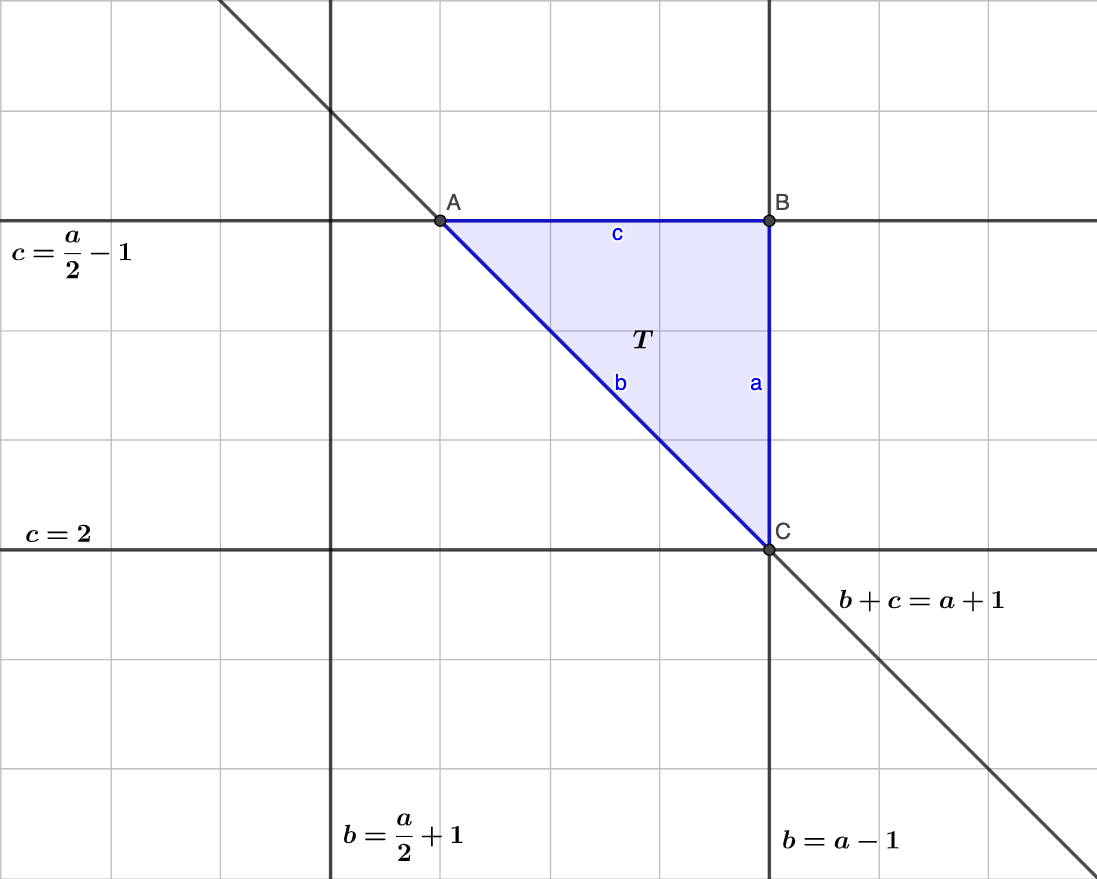}
	\nopagebreak\captionof{figure}{The domain (in blue) of $\ell_1$.}

	\end{center}
where $A=({a\over 2}+2,{a\over 2}-1),B=(a-1,{a\over 2}-1)$ and $C=(a-1,2)$.\\
We have $\nabla\ell_1(b,c)=(a-1-c,a+1-b)$, so that $\nabla\ell_1(b,c)=(0,0)\Leftrightarrow(b,c)=(a+1,a-1)\notin T$. Hence the minimum of $\ell_1$ is reached on the boundary of $T$. We have:
$$\ell_1(b,c)|_{AB}=\ell_1\left( b,{a\over 2}-1\right)={ab+a^2-3a\over 2}\Rightarrow \min_{AB}\ell_1=\ell_1\left( {a\over 2}+2,{a\over 2}-1\right)={3a^2-2a\over 4}$$
$$\ell_1(b,c)|_{BC}=\ell_1(a-1,c)=2c+a^2-3a+2\Rightarrow\min_{BC}\ell_1=\ell_1(a-1,2)=a^2-3a+6$$
$$\ell_1(b,c)|_{CA}=\ell_1(b,a+1-b)=b^2+(-a-3)b+a^2+a+2\Rightarrow\min_{CA}\ell_1=\ell_1\left( {a\over 2}+2,{a\over 2}-1\right) ={3a^2-2a\over4}$$
Under our assumptions it easy to check that $\min\{{3a^2-2a\over 4},a^2-3a+6\}={3a^2-2a\over 4}$ for $a\neq 5$, and since we are assuming $a$ even, we conclude that $\min\ell_1={3a^2-2a\over 4}$. 
\a Now we find the minimum of $\ell_2$. The domain of $\ell_2$ is defined by the following inequalities:
$${a+1\over 2}<b<a-1,\quad c={a\over 2},\quad b+c\geq a+1 \quad ({\rm hence} \;\;   b\geq c)$$
that, under our assumptions, can be refined as follows:
$${a\over 2}+1\leq b\leq a-2,\quad  c={a\over 2}.$$
Hence we get: 
$$\ell_2(b,c)=\ell_2\left( b,{a\over 2}\right) ={(a-2)b+a^2-a\over 2}\Rightarrow\min\ell_2=\ell_2\left( {a\over 2}+1,{a\over 2}\right)={3a^2-2a-4\over 4} .$$
Finally we find the minimum of $\ell_4$. The domain of $\ell_4$ is defined by the following inequalities:
$${a+1\over 2}\leq b<a-1,\quad {a+1\over 2}\leq c<a-1,\quad b+c\geq a+1,  \quad   b\geq c$$
that, under our assumptions, can be refined as follows:
$${a\over 2}+1\leq b\leq a-2,\quad {a\over 2}+1\leq c\leq a-2,\quad b+c\geq a+1,  \quad   b\geq c.$$
The plane region corresponding to these inequalities is the triangle $T'$ as shown by Fig. 3:
	\begin{center}
	\includegraphics[scale=0.7]{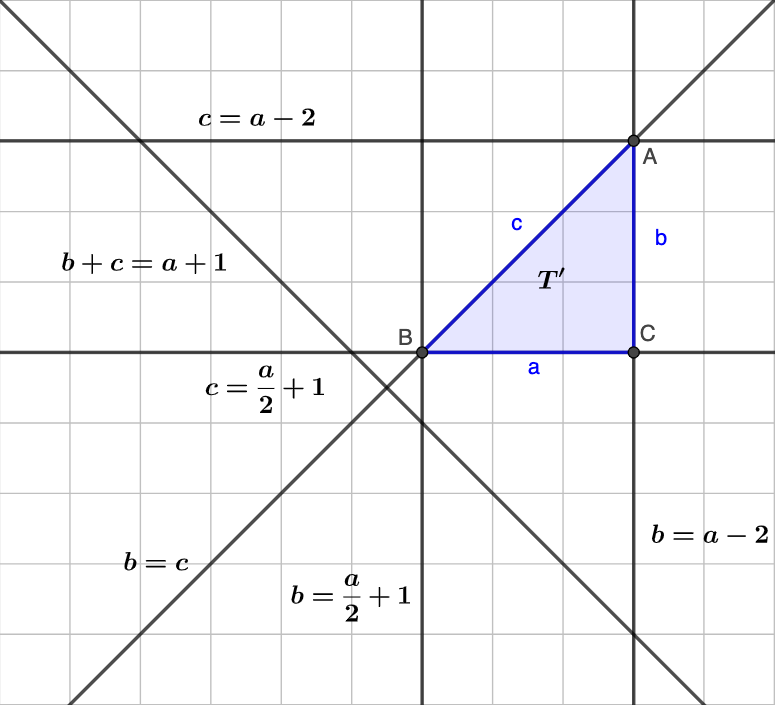}
	\nopagebreak\captionof{figure}{The domain (in blue) of $\ell_4$.}
\end{center}
where $A=(a-2,a-2), B=({a\over 2}+1,{a\over 2}+1)$ and $C=(a-2,{a\over 2}+1)$.\\
Proceeding as we did for $\ell_1$ it easy to check that the minimum of $\ell_4$ is reached on the boundary of $T'$ and it is $\min\ell_4={3a^2-12\over 4}$. 
\dd Finally, an easy computation shows that $$\min\left\lbrace {3a^2-2a\over 4},{3a^2-2a-4\over 4},{3a^2-12\over 4}\right\rbrace={3a^2-2a-4\over 4}$$ 
 and obviously this minimum is attained where $\ell_2$ attains its minimum, i.e. on $ ({a\over 2}+1,{a\over 2})$, 
hence the result is proved for $a$ even. 
\b If $a$ is odd the proof is analogous (but it is easier, since it is not necessary to consider $\ell_2$). \prfend

\begin{remark}\label{holes} \rm Theorem 5.2 of \cite{LP} guarantees that, for an ordinary singularity of multiplicity $a$, all the values in the interval $\left\lfloor {3a^2-2a-4\over 4}\right\rfloor \leq \t  \leq (a-1)^2$ occur. Unfortunately, though
Theorem \ref{poss} allows to recover many such values,  the  curves $C_{b,c} $ do not cover, in general, all the range. Let us see two examples.

\begin{enumerate}[label=(\roman*),leftmargin=*]
	\item $a=9$ corresponds to the interval $55\leq \t \leq 64$. Using Proposition \ref{lunghezze} it is easy to see that the curves $C_{b,c}: x^9+y^9+x^by^c=0$ give, as $b$ and $c$ vary, all these values of $\t$ except for $59$. To find $\t=59$ we can for example consider the curve $D: x^9+y^9+x^5y^7+x^7y^4=0$; a computation with CoCoA shows that  $\t=59$.
		\item $a=10$ corresponds to the interval $69\leq \t \leq 81$ and again it is easy to see that the curves $C_{b,c}: x^{10}+y^{10}+x^by^c=0$ give, as $b$ and $c$ vary, all the values of $\t$ except for $71$ and $74$. To find these two values we can for example consider the curves  $D_1: x^{10}+y^{10}+x^3y^8+x^7y^5=0$ and $D_2: x^{10}+y^{10}+x^2y^9+x^7y^6=0$; $\t$ is 71 for $D_1$ and 74 for  $D_2$.
\end{enumerate}

\end{remark}

\section{Jacobian schemes at double points}

We now consider the Jacobian scheme at a double, not necessarily ordinary, point. Recall that, if $k\geq 1$, a double point of type  $A_{2k-1}$, respectively $A_{2k}$, for a plane curve $C$ is a 2-branched, respectively 1-branched double point which needs  $k$ successive blow ups to be smoothed (hence $A_1$ is a node, $A_2$ an ordinary cusp, $A_3$ a  tacnode and so on). 
 For any $n$, even or odd, the normal form for an $A_n$ singularity is: $y^2-x^{n+1} =0$; in other words, the germ of $C$ at $O$ is analytically equivalent to the germ of the curve $y^2-x^{n+1} =0$ at $O$.
 \b Let $C: f=0$ be a plane curve in $\A^2(\C)$ with a singular point at $O$; 
in \cite{GLS 1}, theorems 2.48 p.150, it is proved that $O$ is double of type $A_n$ if and only if $crk (f)\leq 1$ and the Milnor number  at $O$ is $n$, where $crk(f):= 2- rank\, H(f,O)$ and $H(f,O)$ is the Hessian matrix of $f$ at $O$.

\begin{thm}\label{double}  A point $P$ is a double point of type $A_n$ for a plane curve $D$ if and only if the Jacobian scheme of $D$ at $P$ is a curvilinear scheme of length $n$.
\b Hence, a double point $P$ for a plane curve $D$ is of type $A_n$  if and only if $\tau(D)_P=n$.
 \end{thm}
 \proof  Let us consider the curve $C_n: y^2-x^{n+1} =0$; $C$ has a double point of type $A_n$ at $O$, and no other singularities in the affine plane; its Jacobian ideal is
$$J=(y^2-x^{n+1}, 2y, (n+1)x^{n})= (y,x^{n})$$
hence its Jacobian scheme is a curvilinear scheme of length $n$. The conclusion follows by Theorem \ref{M-Y algebrico}.
 \prfend

\b Theorem \ref{double} gives the following algorithm:

\begin{algorithm}[H]\caption{Study of the double points of a plane curve $C: \sum a_{ij}x_0^ix_1^jx_2^{d-i-j}$}\label{Algorithm1True}
\a  \textbf{Input}: $F =  \sum a_{ij}x_0^ix_1^jx_2^{n-i-j} \in \mathbb{C}[x_0,x_1,x_2]_n$, $n>0$, $P=[a,b,c]$, $F(P)=0$.\\
 \a \textbf{Output}: State if $P$ is a double point for $C$ and its kind: $A_m$.

\a \begin{boxedminipage}{130mm}
    \begin{algorithmic}[1]

	\STATE\label{Alg1Step0} {\it STEP 0)} \  Perform a linear change of coordinates so to have $P = [0,0,1]$; work in the affine chart $\{x_2\neq 0\}$, with $f =  \sum a_{ij}x^iy^j$, and $a_{00}=0$. Set $A:=\C[x,y]$.
	
	   If $(a_{01},a_{10})\neq(0,0)$: $P$ {\bf is a simple point for} $C$, with tangent $a_{01}x+a_{10}y=0$. {\bf STOP}
	
	   If $(a_{01},a_{10})=(0,0)$, and $(a_{20},a_{11},a_{02})=(0,0,0)$: $P$ {\bf is a point of multiplicity $\geq 3$ for} $C$. {\bf STOP}. Otherwise go to the next step
	
	  \a  \STATE\label{Alg1Step1}{\it STEP 1)} \ If $(a_{01},a_{10})=(0,0)$, and $(a_{20},a_{11},a_{02})\neq (0,0,0)$: $P$ {\bf is a double point for} $C$. Go to next step

	 \a \STATE\label{Alg1Step2}{\it STEP 2)} Compute $J = (f, f_x, f_y)$ and  $\dim(A/(J+(x,y)^2)) := \alpha_2$. Go to next step
		
\medskip	
For $r\geq 3$, let 	
 \STATE\label{Alg1Stepr}{\it STEP r)} \ Compute  $\dim(A/J+(x,y)^{r}) := \alpha_{r}$.

  If $\alpha_r = \alpha_{r-1}$, {\bf STOP}: $ \tau$ at $P$ is $\alpha_{r-1}$ and $P$ is of type $A_{r-1}$.

 If $\alpha_r \neq \alpha_{r-1}$ go to next step.

  \end{algorithmic}
\end{boxedminipage}
\end{algorithm}

\section{A remark on the Tjurina number}

In this short section we give a condition for a curve $C$ to have only nodes, through the global Tjurina number $\tau(C)$.

\begin{prop}\label{length} Let $C: f(x,y)=0$ be a plane curve and assume $P\in C$, $m=m_P(C)\geq 2$ and $P$ is not a node; then $X_P$ contains properly the 0-dimensional scheme $(m-1)P$, so that 
$$\ell(X_P)>{m\choose 2}$$
\end{prop} 
\proof First assume that $P$ is an ordinary singularity with $m\geq 3$; then, $X_P \supseteq (m-1)P$ by Lemma \ref{0}, and since for $m\geq 3$ one has ${m\choose 2} < \left\lfloor {3m^2-2m-4\over 4}\right\rfloor$, we get  $X_P \supset (m-1)P$.

\dd If $m\geq 2$ and $P$ is not an ordinary singularity, we can assume that $P=(0,0)$ and the tangent cone contains the double line supported on the $x$-axis, so that 
$$f= y^2h_1\cdot \dots \cdot h_{m-2} + \phi$$
where $h_1,\dots,h_{m-2}$ are linear forms and $\phi$ is sum of forms of degree $\geq m+1$. We have 
$$f_x= \underbrace{y^2(h_1\cdot \dots \cdot h_{m-2})_x}_{\deg m-1} + \phi_x, \qquad f_y= \underbrace{2y(h_1\cdot \dots \cdot h_{m-2})+ y^2(h_1\cdot \dots \cdot h_{m-2})_y }_{\deg m-1}+ \phi_y$$
with $\phi_x=0$ or $\deg \phi_x \geq m$, and $\phi_y=0$ or $\deg \phi_y \geq m$.
Hence the curves $$C_x: f_x=0, \;\;C_y:=f_y=0$$
have a singularity at $P$ with  $m_P(C_x)= m-1$, $m_P(C_y)= m-1$, and the two curves have a common tangent at $P$, i.e. the $x$-axis, so that $C_x \cap C_y$ contains the 0-dimensional scheme $Y$ union of $(m-1)P$ and of the curvilinear scheme of length $m$ supported on the $x$-axis, i.e. $I_Y=(x,y)^{m-1}\cap (x^m,y)=(x^m,x^{m-2}y, \dots,xy^{m-2}, y^{m-1})$; we have $\ell (Y)={m\choose 2}+1$. Moreover, $I_Y \supset (x,y)^m$, that is, $Y \subset mP$.
\dd  The Jacobian scheme at $P$ is the schematic intersection $C \cap C_x \cap C_y$ at $P$, and $C \supset mP \supset Y$, so we get that $Y \subseteq X_P$,  this giving
$\ell(X_P)\geq {m\choose 2}+1$.
\prfend

\begin{prop}\label{rat} Let $C$ be an irreducible curve of degree $d$ and geometric genus $g$ in $\P^2$, with no infinitely near singular points. Then $C$ has only nodes  if and only if $\tau(C)= {d-1 \choose 2}-g$.
\end{prop}
\proof  Let $P_1,\dots,P_r$ be the singular points of $C$, of multiplicity $m_1,\dots,m_r$, and let $\Bbb X$ be its Jacobian scheme. Assume $\tau(C)= {d-1 \choose 2}-g$. Since there are no infinitely near singular points, $g= p_a(C)- \sum_{i=1}^r {m_i\choose 2}$, that is
$${d-1 \choose 2}-g =\sum_{i=1}^r {m_i\choose 2}.$$
Hence the assumption $ {d-1 \choose 2} - g= \tau(C)= \sum_{i=1}^r \ell( \X_{P_i})$ gives:
$$ \sum_{i=1}^r {m_i\choose 2}= \sum_{i=1}^r \ell( \X_{P_i}).$$
If $P_i$ is a node $\ell( \X_{P_i})=1={m_i\choose 2}$, while if $P_i$ is not a node we have $\ell( \X_{P_i})>{m_i\choose 2}$ by Proposition \ref{length}, hence all the singular points must be nodes. The viceversa is immediate. \prfend

\a {\bf Acknowledgements}
\a This study was carried out within the ``0-Dimensional Schemes, Tensor Theory, and Applications" project 2022E2Z4AK – funded by European Union – Next Generation EU  within the PRIN 2022 program (D.D. 104 - 02/02/2022 Ministero dell’Università e della Ricerca).
The first author was partially supported by Politecnico di Torino, the second and third authors were partially supported by Università di Bologna.

\dd {\it Stefano Canino}, Dipartimento di Scienze Matematiche, Politecnico di Torino, Torino, Italy, \linebreak stefano.canino@polito.it; member of INDAM-GNSAGA

\dd {\it Alessandro Gimigliano}, Dipartimento di Matematica, Università di Bologna, Bologna, Italy, \linebreak  alessandr.gimigliano@unibo.it; member of INDAM-GNSAGA

\dd {\it Monica Idà}, Dipartimento di Matematica, Università di Bologna, Bologna, Italy, \linebreak monica.ida@unibo.it

\begin{thebibliography} {biblio} 
\bibitem[AABM]{AABM} M. Alberich-Carrami\~{n}ana, P. Almir\'on, G. Blanco, A. Melle-Hern\'andez, {\it The minimal Tjurina number of irreducible germs of plane curve singularities}. Indiana Univ. Math J. 70, 2001, 3693-3707.
\bibitem[A] {A} P. Almir\'on, {\it On the quotient of Milnor and Tjurina numbers for two-dimensional isolated hypersurface
singularities}. Mathematische Nachrichten 295, 2022, 1254–1263. DOI: 10.1002/mana.202100371.
\bibitem[AM] {AM}  M. F. Atiyah, I. G. Macdonald, {\it Introduction to Commutative Algebra}. Addison-Wesley Publishing Company, Reading 1969.
\bibitem[BGM] {BGM} J. Brian\c{c}on, M. Granger and P. Maisonobe, {\it Le nombre de modules du germe de courbe plane $x^a + y^b = 0$}. Math. Ann. 279, 535–551 (1988).
\bibitem[BK] {BK} E. Brieskorn, H. Knörrer, {\it  Plane algebraic curves}. Birkhäuser Verlag, Basel, 1986
\bibitem[CLO] {CLO}  D. A. Cox, J. Little, D. O'Shea {\it Ideals, varieties and Algorithms}. Springer 2015.
\bibitem[CoCoa] {CoCoa} J. Abbott, A. M. Bigatti, L. Robbiano {\it CoCoA: a system for doing Computations in Commutative Algebra.}
Available at http://cocoa.dima.unige.it
\bibitem[Fu]{F} W. Fulton, {\it Algebraic Curves}. Benjamin/Cummings Publishing Company, Reading, Massachusetts 1969.
\bibitem[GH]{GH} Y. Genzmer, M. E. Hernandes {\it On the Saito basis and the Tjurina number for plane branches.} Transactions of the AMS, 2020, 373 (5), 3693-3707. 10.1090/tran/8019. hal-03794511 
\bibitem[GI] {GI} A. Gimigliano, M. Idà, {\it Remarks on double points of plane curves}. Geometriae Dedicata 217, 2023.
\bibitem[G]{G}{J. M. Granger}, {\it Sur une espace de modules de germe de courbe plane.} Bull. Sci. Math. 26 s6r. 103,1979.
\bibitem[GLS 1]{GLS 1} G.-M.  Greuel, C. Lossen, E. Shustin, {\it Introduction to Singularities and Deformations}. Springer 2007.
\bibitem[GLS 2]{GLS 2} G.-M.  Greuel, C. Lossen, E. Shustin, {\it Singular Algebraic Curves} With an Appendix by Oleg Viro. Springer 2018.
\bibitem[Ha]{Ha}  R. Hartshorne, {\it Algebraic Geometry}. Springer, Grad. Texts in Math. 52, Berlin, Heidelberg, New York 1977.
\bibitem[HH]{HH}  A. Hefez, M. E. Hernandes, {\it Analytic classification of plane branches up to multiplicity} 4. J. of Symbolic Computation 44, 2009, 626-634.
\bibitem[LP]{LP} O. A. Laudal, G. Pfister, {\it Local moduli and singularities}. Springer, Lecture notes in Mathematics 1310, New York/Berlin 1988.
\bibitem[MY]{MY} J. N. Mather, S. S.-T. Yau, {\it Classification of Isolated Hypersurface Singularities by Their Moduli Algebras}. Invent. Math. 69, 1982, 243-251.
\bibitem[GAGA] {GAGA} J.-P. Serre, {\it Géométrie algébrique et géometrie analytique}. Annales de l'institut Fourier, tome 6, 1956, p.1-42
\bibitem[W]{W}{M. Watari}, {\it Plane curve singularities whose Milnor and Tjurina numbers differ by three}. Advanced Studies in Pure Mathematics 46, 2007
Singularities in Geometry and Topology 2004, 273-298.

\end{thebibliography}
\end{document}